\documentclass[aop,MSNbibl,secthm,seceqn,citesort,dvips]{arximspdf}

\doi{10.1214/10-AOP604}
\volume{39}
\issue{6}
\pubyear{2011}
\firstpage{2474}
\lastpage{2496}

\makeatletter
\newtheorem{prop}{Proposition}[section]
\newtheorem{cor}{Corollary}[section]

\newcommand{\afrac}[2]{{(#1/#2)}}
\newcommand{\bfrac}[2]{{(#1)/#2}}
\newcommand{\cfrac}[2]{{#1/#2}}
\newcommand{\widehata}[1]{\hat#1}
\newcommand{\widehatb}[1]{\hat#1}
\newcommand{\vertvertvert}{\|\!|}

\newcommand{\R}{\mathbb{R}}
\newcommand{\nN}{n \in\mathbb{N}}
\newcommand{\C}{\mathbb{C}}
\newcommand{\Tr}{\operatorname{Tr}}
\newcommand{\Dom}{\operatorname{Dom}}
\newcommand{\G}{\widehat{G}}
\newcommand{\CID}{\operatorname{CID}_{\R}(G)}
\newcommand{\hs}{\mathrm{hs}}
\newcommand{\cal}{\mathcal}
\newproclaim{examplea}{Example}
\newproclaim{exampleb}{Example}
\newproclaim{example}{Example}
\newtheorem{theorem}{Theorem}[section]
\newproclaim{pff}{Proof}
\makeatother

\begin{document}
\begin{frontmatter}

\title{Infinitely divisible central probability measures on compact Lie groups---regularity,
semigroups and transition kernels}
\runtitle{Infinitely divisible central probability measures}

\begin{aug}
\author[A]{\fnms{David} \snm{Applebaum}\corref{}\ead[label=e1]{D.Applebaum@sheffield.ac.uk}}
\runauthor{D. Applebaum}
\affiliation{University of Sheffield}
\address[A]{Department of Probability\\ \quad and Statistics\\
University of Sheffield\\
Hicks Building, Hounsfield Road\\
Sheffield S3 7RH\\
United Kingdom\\
\printead{e1}} 
\end{aug}

\received{\smonth{7} \syear{2010}}
\revised{\smonth{8} \syear{2010}}

%
\begin{abstract}
We introduce a class of central symmetric infinitely divisible
probability measures on compact Lie groups by lifting the
characteristic exponent from the real line via the Casimir
operator. The class includes Gauss, Laplace and stable-type
measures. We find conditions for such a measure to have a smooth
density and give examples. The Hunt semigroup and generator of
convolution semigroups of measures are represented as
pseudo-differential operators. For sufficiently regular
convolution semigroups, the transition kernel has a~tractable
Fourier expansion and the density at the neutral element may be
expressed as the trace of the Hunt semigroup. We compute the short
time asymptotics of the density at the neutral element for the
Cauchy distribution on the $d$-torus, on SU(2) and on SO(3), where
we find markedly different behaviour than is the case for the
usual heat kernel.
\end{abstract}

%

\end{frontmatter}

\section{Introduction}\label{sec1}

The heat kernel on a compact Riemannian manifold has been the subject
of extensive
investigations by analysts, geometers and probabilists. One reason for
this is that its
small and large time asymptotic expansions contain important
topological and geometric
information (see, e.g., \cite{Ros}). Another reason is that it is the
transition density of
manifold-valued Brownian motion which is a stochastic process of
intrinsic interest
(see, e.g., \cite{elw}). If the manifold is a Lie group, then the heat
kernel is naturally related to Dedekind's eta function via Macdonald's
identities (see \cite{Feg1}). In this paper, we will mainly be
concerned with compact Lie groups. In this case, the heat kernel is
naturally associated to a vaguely (or equivalently, weakly) continuous
convolution semigroup of probability measures which we'll refer to as
the ``heat semigroup'' in the sequel.

The study of the entire class of such convolution semigroups has
had a~long development (see, e.g., \cite{Hu,He1,App2}). From a probabilistic point of view, they correspond
to L\'{e}vy processes, that is, stochastic processes with stationary
and independent increments. Compared to Brownian motion which has
continuous sample paths (with probability one), the paths of the
generic L\'{e}vy process are only right continuous and have jump
discontinuities of arbitrary size occurring at random times.

The purpose of this paper is to study a class of convolution
semigroups which on the one hand, are sufficiently close in
structure to the heat semigroup to enable us to do some
interesting analysis and on the other hand, are sufficiently broad
as to display all the interesting features that one finds with a
generic L\'{e}vy process. The first observation is that measures
comprising the heat semigroup are central and so we focus on this
class. It is worth pointing out that that central measures as a
class have also received some attention from analysts (see, e.g.,
\cite{Rag,Hare}). Second, we remark that if $(\mu_{t}, t
\geq0)$ is the heat semigroup then its noncommutative Fourier
transform (see \cite{He2,Sieb} for background on this
concept) takes the form $\widehata{\mu_{t}}(\pi) =
e^{-\afrac{t}{2}\kappa_{\pi}}I_{\pi}$ at the irreducible
representation $\pi$ where $-\kappa_{\pi}I_{\pi}$ is the Casimir
operator. But $u \rightarrow u^{2}/2$ is the negative-definite
function on the real line associated with the standard Gaussian
measure. The generalization that we make here is to consider a~%
class of semigroups that are given by the prescription
$\widehata{\mu_{t}}(\pi) =
e^{-t\eta(\kappa_{\pi}^{\cfrac{1}{2}})}I_{\pi}$ where $\eta$ is a
symmetric negative definite function. Other examples of measures
subsumed within this class include the Laplace distribution on a
Lie group, which has been untilized in recent statistical work on
the problem of deconvolution (see, e.g., \cite{KK,KR}) and
analogues of stable laws. Indeed any semigroup of probability
measures that is obtained by subordinating the heat semigroup
belongs to this class. We study these measures both from the
static perspective, where the emphasis is on a single infinitely
divisible measure, and the dynamic perspective where we focus on
an entire semigroup.

The organization of this paper is as follows. In Section~\ref{sec2}, we study
central probability
measures, introduce our main class and examine some examples. In
Section~\ref{sec3}, we use
Sobolev spaces to find conditions on our induced measures which enable
them to have a
smooth density. In Section~\ref{sec4}, we turn our attention to convolution
semigroups and the
associated semigroup of operators (the \textit{Hunt semigroup}) on the
$L^{2}$ space of
normalized Haar measure. When $G$ is a Euclidean space, it is known
(see \cite{Jac1}, Chapter 3 of~\cite{Abook}) that these operators, and their generators, can be
realized as
pseudo-differential operators. Using Peter--Weyl theory, Ruzhansky and
Turunen \cite{RT}
have developed an intrinsic theory of pseudo-differential operators on
compact groups. We adapt this
theory to our needs and show that the Hunt semigroup and its generator are
pseudo-differential operators in the sense of Ruzhansky and
Turunen.\vadjust{\goodbreak}
This part of the work
is carried out in full generality. In the case of our induced class, we
show that the
generator has the same Sobolev regularity as the Laplacian. Although
these results have an
analytic flavor, they are important for probabilists as they indicate a
route to
investigate general classes of Feller--Markov processes on compact Lie
groups using the
symbol of the generator as the key tool (see \cite{Jac4} for an
account of this theory in
the case where $G$ is Euclidean space).

In Section~\ref{sec5}, we investigate the transition kernel for convolution
semigroups of central measures. We remark that the first
investigation of densities for such measures (under a
hypo-ellipticity condition that we do not require here) were made
by Liao in \cite{Liao1} (see also Theorem~4.4 in \cite{Liao},
page~96). A~necessary and sufficient condition for the semigroup to
be trace-class for any positive time is that the corresponding
probability measure has a square-integrable density \cite{App2}.
We compute the trace in both the main $L^{2}$-space and the
subspace of central functions. In the former case, the coordinate
functions form a complete set of eigenvectors for the Hunt
semigroup. Comparing the traces in these two spaces, leads to an
interesting inequality for transition kernels which appears to be
new even in the heat kernel case. Finally, in Section~\ref{sec6}, we study
the small time asymptotics of the transition kernel in the case of
the Cauchy distribution on the $d$-torus, on $\operatorname{SU}(2)$ and on
$\operatorname{SO}(3)$ and show that it blows up much faster than the heat
kernel.

\vspace*{-3pt}\section{Infinite divisibility of central measures}\label{sec2}

Let $G$ be a compact group with neutral element $e$ and let ${\cal
M}(G)$ be the set of all probability measures defined on $(G,
{\cal B}(G))$ where ${\cal B}(G)$ is the Borel $\sigma$-algebra of
$G$. We say that $\mu\in{\cal M}(G)$ is \textit{central} (or \textit{conjugate-invariant}) if $\mu(\sigma A \sigma^{-1}) = \mu(A)$ for
all $\sigma\in G, A \in{\cal B}(G)$ and $\mu$ is said to be \textit{symmetric} if $\mu(A^{-1}) = \mu(A)$ for all $A \in{\cal B}(G)$.
Let ${\cal M}_{c}(G)$ (${\cal M}_{s}(G)$) be the subsets of ${\cal
M}(G)$ comprising central (symmetric) measures (resp.) and
define ${\cal M}_{c,s}(G):= {\cal M}_{c}(G)\cap{\cal M}_{s}(G)$.
Normalized Haar measure on $G$ will always be denoted $d\sigma$
when integrating functions of $\sigma\in G$.

Let $\widehat{G}$ be the set of all equivalence classes of
irreducible representations of~$G$. We will, without further comment,
frequently identify equivalence classes
with a particular representative element when there is no loss of
generality. The trivial representation
will always be denoted by $\delta$. Each $\pi\in\widehat{G}$
acts as a $d_{\pi} \times d_{\pi}$ unitary matrix on a complex linear
space $V_{\pi}$ having dimension~$d_{\pi}$. We define the Fourier
transform of each $\mu\in
{\cal M}(G)$ to be the Bochner integral
\[
\widehata{\mu}(\pi) = \int_{G}\pi(\sigma)\mu(d\sigma),
\]
where $\pi\in\widehat{G}$. We will frequently use the well-known and
easily verified fact that
\[
\widehat{\mu* \nu}(\pi) = \widehata{\mu}(\pi)\hat{\nu}(\pi)
\]
for all $\mu, \nu\in{\cal M}(G), \pi\in\G$, where $*$ denotes
convolution of\vadjust{\goodbreak} measures.

Suppose we are given $\mu\in{\cal M}(G)$. It is shown in
\cite{SLLM} that $\mu\in{\cal M}_{c}(G)$ if and only if for each
$\pi\in\widehat{G}$ there exists $c_{\pi} \in\C$ such that
\begin{equation} \label{conj1}
\widehata{\mu}(\pi) = c_{\pi}I_{\pi},
\end{equation}
where $I_{\pi}$ is the identity matrix acting on $V_{\pi}$.
Indeed this is a straightforward consequence of Schur's lemma.
Moreover, one has the formula
\begin{equation}
c_\pi= \frac{1}{d_\pi} \int_{G} \chi_{\pi}(\sigma)\mu(d\sigma),
\end{equation}
where $\chi_{\pi}(\cdot): = \operatorname{tr}(\pi(\cdot))$ is the group
character.

It is well known (and easily verified) that $\mu\in{\cal
M}_{s}(G)$ if and only of $\widehata{\mu}(\pi)$ is self-adjoint for
all $\pi\in\widehat{G}$. Consequently, $\mu\in{\cal M}_{cs}(G)$
if and only if $\widehata{\mu}(\pi) = c_{\pi}I_{\pi}$ with $c_{\pi
} \in\R$ for all $\pi\in\G$.

A probability measure $\mu$ is \textit{infinitely divisible} if for
each $\nN$ there exists $\nu_{n} \in{\cal M}(G)$ such that
$\nu_{n}^{*n} = \mu$. In this case, we write $\mu_{\cfrac{1}{n}} : =
\nu_{n}$.


\begin{prop} \label{non-zero}
If $G$ is a compact Lie group and $\mu\in
M_{cs}(G)$ is infinitely divisible, then
for each $\pi\in\widehat{G}$ there exists $\alpha_{\pi} \leq0$
such that $\widehata{\mu}(\pi) = e^{\alpha_{\pi}}I_{\pi}$.
\end{prop}

\begin{pf}  By the results on pages~220--221 of \cite{He1}, $\mu$ may
be embedded as~$\mu_{1}$ into a
vaguely continuous convolution semigroup of probability measures
$(\mu_{t}, t \geq0)$ where $\mu_{0}$ is normalized Haar measure on a
closed subgroup~$H$ of~$G$. It follows
(see \cite{App1,LoNg}) that for each $\pi\in\G$,
$(\widehata{\mu_{t}}(\pi), t
\geq0)$ is a~strongly continuous contraction semigroup of
matrices acting on $V_{\pi}$ and so we may write
$\widehata{\mu_{t}}(\pi) = \widehata{\mu_{0}}(\pi)e^{t A_{\pi}}$
for all $t \geq0$ where
$A_{\pi}$ is a $d_{\pi} \times d_{\pi}$ matrix. Now since $\mu_{1}
\in M_{cs}(G)$, there exists $\lambda_{\pi} \in\R$ such that
\renewcommand{\theequation}{*}
\begin{equation}\label{eq*}
\widehata{\mu_{1}}(\pi) = \widehata{\mu_{0}}(\pi)e^{A_{\pi}} =
\lambda_{\pi}I_{\pi}  \ldots.
\end{equation}
\renewcommand{\theequation}{\arabic{section}.\arabic{equation}}
\setcounter{equation}{2}
If $\lambda_{\pi} = 0$, the required result holds with $\alpha_{\pi
} = -\infty$ so assume that $\lambda_{\pi} \neq0$.
Since $\mu_{1} = \mu_{1} * \mu_{0}$, we have
\[
\widehata{\mu_{0}}(\pi)e^{A_{\pi}}\widehata{\mu_{0}}(\pi) =
\lambda_{\pi}I_{\pi}.
\]
On the other hand, post-multiplying both sides of (\ref{eq*}) by $\widehata{\mu
_{0}}(\pi)$ yields
\[
\widehata{\mu_{0}}(\pi)e^{A_{\pi}}\widehata{\mu_{0}}(\pi) =
\lambda_{\pi}\widehata{\mu_{0}}(\pi).
\]
It follows that $\widehata{\mu_{0}}(\pi) = I_{\pi}$ and hence $H = \{
e\}$. We then have
$A_{\pi} = \alpha_{\pi}I_{\pi}$ where $\alpha_{\pi} \in\R$ and
$\lambda_{\pi} = e^{\alpha_{\pi}}$. But
$\widehata{\mu_{1}}(\pi)$ is a contraction on $V_{\pi}$ and hence
$\alpha_{\pi} \leq0$.
\end{pf}





\begin{examplea*}[(The compound Poisson distribution)] Consider
the probability measure $\mu_{\lambda, \gamma}$ where $\gamma$ is
a given probability measure on $G$ and $\lambda> 0$. This is defined by
\[
\mu_{\lambda, \gamma}: =
e^{-\lambda}\sum_{n=0}^{\infty}\frac{\lambda^{n}}{n!}\gamma^{*n}.\vadjust{\goodbreak}
\]

It is well known (see, e.g., \cite{SLLM}) that for all $\pi\in
\widehat{G}$,
\[
\widehata{\mu_{\lambda, \gamma}}(\pi) = \exp{\bigl\{\lambda\bigl(\widehatb{\gamma}(\pi) -
I_{\pi}\bigr)\bigr\}}.
\]
\end{examplea*}

\begin{prop}
\begin{longlist}[(2)]
\item[(1)] The measure $\mu_{\lambda, \gamma}$ is central if and only if
$\gamma$
is.
\item[(2)] The measure $\mu_{\lambda, \gamma}$ is symmetric if and
only if $\gamma$ is.
\end{longlist}
\end{prop}

\begin{pf}
\begin{longlist}[(2)]
\item[(1)] The if part is
straightforward and is established in Proposition~4 of~%
\cite{SLLM}. Conversely, if $\mu_{\lambda, \gamma}$ is
central then for all $g \in G, \pi\in\widehat{G}$
\[
\pi(g) \widehata{\mu_{\lambda, \gamma}}(\pi)\pi(g^{-1}) =
\widehata{\mu_{\lambda, \gamma}}(\pi),
\]
and so
\[
\widehata{\mu_{\lambda, \gamma}}(\pi) = \exp{\bigl\{\lambda\bigl(\pi
(g)\widehatb{\gamma}(\pi)\pi(g^{-1}) -
I_{\pi}\bigr)\bigr\}}.
\]
Now by uniqueness of Fourier transforms and injectivity of the
exponential map on matrices, we have
\[
\pi(g)\widehatb{\gamma}(\pi)\pi(g^{-1}) =
\widehatb{\gamma}(\pi)
\]
for all $g \in G, \pi\in\widehat{G}$ and the result follows.
\item[(2)] This result is proved similarly using the fact that a
probability measure is symmetric if and only if its Fourier
transform comprises self-adjoint matrices.
\end{longlist}
\upqed
\end{pf}

It follows that a central probability measure $\mu$ is a compound
Poisson distribution if and only if there exists $\lambda> 0$ and a
central probability measure~$\gamma$ with $\widehatb{\gamma}(\pi) =
b_{\pi}I_{\pi}$ for all $\pi\in\widehat{G}$ such that
\begin{equation} \label{cpd}
\widehata{\mu}(\pi) = \exp\{\lambda(b_{\pi} - 1) I_{\pi}\}.
\end{equation}

We now introduce a class of central symmetric measures which are
key to this paper. For this part, we assume that $G$ is a compact
Lie group. Let $\rho$ be a symmetric infinitely divisible
probability measure on $\R$. Then we have the L\'{e}vy--Khintchine
formula
\[
\int_{\R}e^{iux}\rho(dx) = e^{-\eta(u)},
\]
where
\begin{equation} \label{Ls}
\eta(u) = \frac{1}{2}\sigma^{2}u^{2} +\int_{\R-\{0\}}\bigl(1-\cos
(uy)\bigr)\nu(dy),
\end{equation}
where $\sigma\geq0$ and $\nu$ is a symmetric L\'{e}vy
measure on $\R-\{0\}$, that is, a~$\sigma$-finite Borel symmetric measure
for which $\int_{\R-\{0\}}\min\{1, |x|^{2}\}\nu(dx) < \infty$
(see, e.g., \cite{BF}).\vadjust{\goodbreak}

For each $\pi\in\widehat{G}$, let $K_{\pi}$ be the Casimir
operator acting in $V_{\pi}$. Then $K_{\pi} = -
\kappa_{\pi}I_{\pi}$ where $\kappa_{\pi} \geq0$ with $\kappa
_{\pi} = 0$ if and only if $\pi= \delta$. If $\mu$ is a
probability measure on $G$ for which
\begin{equation} \label{mainc}
\widehata{\mu}(\pi) = e^{-\eta(\kappa_{\pi}^{\cfrac{1}{2}})}I_{\pi},
\end{equation}
we say that $\mu$ is a central symmetric probability measure on
$G$ \textit{induced by an infinitely divisible probability measure}
on $\R$ and we write $\mu\in \operatorname{CID}_{\R}(G)$.\

The following two examples have been applied to statistical inference
on groups (see, e.g., \cite{KR,KK}).

\begin{example}[(Gaussian measure)] Here we take $\nu=0$ and so
$c_{\pi} =\break \exp \{-\frac{1}{2}\times \sigma^{2}\kappa_{\pi} \}$.
Gaussian measure is embeddable into the Brownian motion or heat
semigroup of measures for which $\widehata{\mu}_{t}(\pi) =
\exp \{-\frac{t}{2}\sigma^{2}\kappa_{\pi} \}$ for $t
\geq
0$ which has been extensively studied by both analysts and
probabilists.
\end{example}

\begin{example}[(The Laplace distribution on $G$)] $\!\!$Here we take $\sigma\,{=}\,0,
\nu(dx)\,{=}\break
\frac{\exp \{-\cfrac{|x|}{\beta} \}}{|x|}\,dx$ (with
$\beta
> 0$) and $c_{\pi} = (1 + \beta^{2} \kappa_{\pi})^{-1}$ (see
\cite{Sa}, page~98 for a~discussion of the underlying distribution on
$\R$).
\end{example}

Now consider a central symmetric compound Poisson
distribution $\mu_{\lambda, \gamma}$. We consider conditions under
which $\mu_{\lambda, \gamma} \in\CID$. First, take $\sigma= 0$
and~$\nu$ to be a finite symmetric measure in (\ref{Ls}) and
rewrite
\[
\eta(u) = \lambda\int_{\R-\{0\}}\bigl(1 -
\cos(uy)\bigr)\tilde{\nu}(dy),
\]
where $\lambda: = \nu(\R-\{0\})$ and
$\tilde{\nu}(\cdot):=\frac{1}{\lambda}\nu(\cdot)$. For
$\mu_{\lambda, \gamma} \in\CID$ with this value of $\lambda$, we
require that $b_{\pi} = g(\kappa_{\pi}^{\cfrac{1}{2}})$ in
(\ref{cpd}) where $g(u) = \int_{\R}\cos(ux)\tilde{\nu}(dx)$. For
example, if we take $\nu$ to be
a constant multiple of a centred Gaussian measure with variance $\sigma
^{2}$ on $\R$ then
$b_{\pi} = \exp \{-\frac{1}{2}\sigma^{2}\kappa_{\pi} \}$.

We now consider an important subclass of measures in
$ \operatorname{CID}_{\R}(G)$. Let $(\rho_{t}^{f}, t \geq0)$ be the law of a
subordinator with associated Bernstein function $f\dvtx (0, \infty)
\rightarrow[0, \infty)$ so that $(\rho_{t}^{f}, t \geq0)$ is a
vaguely continuous convolution semigroup of probability measures
on $[0, \infty)$ and for each $t \geq0, u > 0$,
\begin{equation} \label{sub}
\int_{0}^{\infty}e^{-us}\rho_{t}^{f}(ds) = e^{-tf(u)},
\end{equation}
and $f$ has the generic form
\[
f(u) = au + \int_{(0, \infty)}(1 - e^{-uy})\lambda(dy),
\]
where $a \geq0$ and $\int_{(0, \infty)}\min\{1, y\}\lambda(dy) <
\infty$ (see, e.g.,~\cite{Sa}, Section~30 and~\cite{Abook}, Section~1.3.2 for details).
 It is straightforward to verify that if
\mbox{$(\mu_{t}, t \geq0)$} is a~vaguely\vadjust{\goodbreak} continuous convolution
semigroups of measures on $G$ and \mbox{$(\rho_{t}^{f}, t \geq0)$} is a
subordinator as above then we get another vaguely continuous
convolution semigroups of measures on $G$ which we denote
$(\mu_{t}^{f}, t \geq0)$ via the vague integral
\[
\mu_{t}^{f}(A) = \int_{0}^{\infty} \mu_{s}(A)\rho_{t}^{f}(ds)
\]
for $A \in{\cal B}(G)$. Now let $(\mu_{t}, t \geq0)$ be the
Brownian motion semigroup with $\sigma= \sqrt{2}$. Then for each
$\pi\in\G, t \geq0$ we have
\begin{eqnarray*}\widehata{\mu_{t}^{f}}(\pi) &
= &
\int_{0}^{\infty}\int_{G}\pi(\sigma)\mu_{s}(d\sigma)\rho
_{t}^{f}(ds)\\[-2pt]
& = &  \biggl(\int_{0}^{\infty}e^{-s\kappa_{\pi}}\rho
_{t}^{f}(ds) \biggr)I_{\pi}\\[-2pt]
& = & e^{-tf(\kappa_{\pi})}I_{\pi},
\end{eqnarray*}
and so $\mu_{1}^{f} \in
\operatorname{CID}_{\R}(G)$ with $\eta(\kappa_{\pi}^{\cfrac{1}{2}}) =
f(\kappa_{\pi})$.

Note that the Laplace distribution (as described above) is
obtained in this way with $f(u) = \log(1 + \beta^{2}u)$. It is
worth pointing out that it also arises as
$\beta^{-2}V^{\beta^{-2}}$ where for $c > 0, V^{c}$ is the
potential measure of the Brownian motion semigroup defined by the
vague integral $V^{c}(\cdot) =
\int_{0}^{\infty}e^{-ct}\mu_{t}(\cdot)\,dt$ (see~\cite{Sa}, pages~203--205 for the case in $\R^{d}$).

Other examples of measures in $\operatorname{CID}_{\R}(G)$ which are
obtained by subordination include stable-type distributions where
$\sigma= 0$ in (\ref{Ls}) and $\nu(dx) =
\frac{b^{\alpha}}{|x|^{1+\alpha}}\,dx$ where $b > 0$ and $0 <
\alpha< 2$. In this case, we have $f(u) =
b^{\alpha}u^{\cfrac{\alpha}{2}}$ and $c_{\pi} =
\exp \{-b^{\alpha}\kappa_{\pi}^{\cfrac{\alpha}{2}} \}$. We
may also consider the relativistic Schr\"{o}dinger distribution
for $m > 0$ where $f(u) = \sqrt{u + m^{2} - m}$ and $c_{\pi} =
e^{-(\sqrt{m^{2} + \kappa_{\pi}} - m)}$. It again has $\sigma= 0$
in (\ref{Ls}). The precise form of $\nu$ is complicated and as we
do not require it here we refer the reader to \cite{Ic}.

It is an interesting problem to determine the class of all $\eta$ in
(\ref{Ls}) which give rise to a probability measure on $G$ of the form
(\ref{mainc}).

\vspace*{-3pt}\section{Regularity of densities}\label{sec3}

In this section we will assume that $G$ is a~compact semi-simple
Lie group having Lie algebra $\mathbf{g}$. We say that $\mu\in{\cal
M}(G)$ has a density $k \in L^{1}(G, \R)$ if $\mu$ is absolutely
continuous with respect to normalized Haar measure on $G$. We then
define $k$ to be the Radon--Nikod\'{y}m derivative
$\frac{d\mu}{d\sigma}$.

If a density $k$ exists for $\mu\in{\cal M}_{c}(G)$ with
$\widehata{\mu}(\pi) = c_{\pi}I_{\pi}$ and $k \in L^{2}(G, \R)$
then it has the form
\begin{equation} \label{dens}
k(\sigma) = \sum_{\pi\in\G}d_{\pi}\overline{c_{\pi}}\chi_{\pi
}(\sigma)
\end{equation}
for almost all $\sigma\in G$ (see equation (3.4) in \cite{App}).

Before we investigate densities in greater detail, we need some preliminaries.\vadjust{\goodbreak}


\subsection{Dominant weights}\label{sec3.1}

Fix a maximal torus $\mathbb{T}$ in $G$. Let $\mathrm{T}$ be its
Lie algebra and $\mathrm{T}^{*}$ be the dual vector space to
$\mathrm{T}$. Let $P$ be the lattice of weights in
$\mathrm{T}^{*}$ and $D \subset\mathrm{T}^{*}$ be the dominant
chamber. The celebrated theorem of the highest weight asserts that
there is a one-to-one correspondence between elements of
$\widehat{G}$ and the highest weights which are precisely the
members of $P \cap D$. For details, see, for example, Chapters~8 and~9 in
\cite{Feg}. In the following, the inner product $\langle\cdot, \cdot
\rangle$ and norm $|\cdot|$ on $\mathrm{T}^{*}$ are that induced by
the Killing form via duality.

Let $\lambda_{\pi}$ be the dominant weight for the representation
$\pi$.
Then we know from Sugira \cite{Sug} [page 39, equation (1.17)] that
\begin{equation} \label{w1}
d_\pi\leq N |\lambda_{\pi} |^m,
\end{equation}
where $N$ is a universal constant and
\begin{equation} \label{rank}
m = \tfrac{1}{2}\bigl(\operatorname{dim}(G) - r\bigr),
\end{equation}
where $r$ is the rank of $G$, that is, the dimension of any
maximal torus. It is also well known that (see, e.g., \cite{Sug},
Lemma~1.1 or \cite{Knapp}, Proposition~5.28)
\begin{equation} \label{eig}
\kappa_{\pi} = |\lambda_{\pi} + \rho|^{2} - |\rho|^{2} = \langle
\lambda_\pi, \lambda_\pi+ 2\rho\rangle,
\end{equation}
where $\rho$ is half the sum of positive roots. It follows easily
that
\begin{equation} \label{w2}
|\lambda_\pi|^2 \leq\kappa_{\pi} \leq|\lambda_\pi|^2 +
2|\lambda_\pi|  |\rho| \leq C(1+|\lambda_{\pi}|^2),
\end{equation}
where $C > 1$ is a constant.

We also need the fact (which is implicit in the proof of Lemma~1.3 in~%
\cite{Sug}) that there exists $C_{1}, C_{2} > 0$ such that for all
$\lambda\in P \cap D$, there exists $n = (n_{1}, \ldots, n_{r}) \in
\mathbb{Z}^{r}$ such that
\begin{equation} \label{inteq}
C_{1}\|n\| \leq|\lambda| \leq C_{2}\|n\|,
\end{equation}
where $\|n\|^{2}:=n_{1}^{2} + \cdots+ n_{r}^{2}$.

The final result we need from Sugiura \cite{Sug} is Lemma~1.3
therein that
\begin{equation} \label{w3}
\zeta(s) : = \sum_{\lambda\in P \cap D - \{0\}} \langle\lambda,
\lambda\rangle^{-s}
\end{equation}
converges if $2s > r$.



\subsection{Sobolev spaces}\label{sec3.2}

Let $\{X_{1}, \ldots, X_{d}\}$ be a basis for the Lie algebra $\mathbf{g}$
of left-invariant vector fields. We define the Sobolev space
${\cal H}_{p}(G)$ by the prescription
%
\[
{\cal H}_{p}(G):=\{f \in L^{2}(G); X_{i_{1}}\cdots X_{i_{k}}f
\in L^{2}(G); 1 \leq k \leq p, i_{1}, \ldots, i_{k} = 1, \ldots,
d\}.
\]
It is a complex separable Hilbert space with associated
norm
\[
\vertvertvert f\vertvertvert _{p}^{2} = \|f\|^{2} + \sum_{i_{1}, \ldots,
i_{k}}\|X_{i_{1}}\cdots
X_{i_{k}}f\|^{2}.
\]

It is not difficult to show that an equivalent norm is given by
\begin{equation} \label{Sob}
\vertvertvert f\vertvertvert _{p}^{2} = \sum_{\pi\in\widehat{G}}d_{\pi}(1 +
\kappa_{\pi})^{p}\operatorname{tr}(\widehatb{f}(\pi)\widehatb{f}(\pi)^{*}),
\end{equation}
where $\widehatb{f}(\pi):=\int_{G}\pi(\sigma^{-1})f(\sigma)\,d\sigma$
is the Fourier transform\footnote{Note that we are here using the
analyst's convention for Fourier transforms of functions which, as
usual, is not quite consistent with the probabilist's convention
for Fourier transforms of measures.} (and we are abusing notation
by using $\vertvertvert \cdot\vertvertvert $ in each case).

As is pointed out in \cite{RT}, Section~10.3.1, ${\cal H}_{p}(G)$
coincides with the usual Sobolev space on a manifold constructed
using partitions of unity. So in particular, the Sobolev embedding
theorem extends to this context and hence
\[
C^{\infty}(G) \supseteq \bigcap_{k \in\mathbb{N}}{\cal H}_{k}(G).
\]

\subsection{A regularity result}\label{sec3.3}

We summarize the results we need on regularity in the following.
\begin{prop} \label{main}
Let $\mu\in{\cal M}_{c}(G)$ with $\widehata{\mu}(\pi) =
c_{\pi}I_{\pi}$ for all $\pi\in\widehat{G}$.
\begin{longlist}[(3)]
\item[(1)] The measure $\mu$ has a square-integrable density if and
only if
\begin{equation} \label{L2}
\sum_{\pi\in\widehat{G}}d_{\pi}^{2}|c_{\pi}|^{2} < \infty.
\end{equation}
\item[(2)] The measure $\mu$ has a continuous density if
\begin{equation}
\sum_{\pi\in\widehat{G}} d^{2}_{\pi} |c_{\pi}| < \infty.
\label{eqn:contdens1}
\end{equation}
\item[(3)] The measure $\mu$ has a $C^{k}$ density if
\begin{equation} \label{Sobcen}
\sum_{\pi\in\widehat{G}}d_{\pi}^{2}(1 +
\kappa_{\pi})^{p}|c_{\pi}|^{2} < \infty,
\end{equation}
where $p > k + \frac{d}{2}.$
\end{longlist}
\end{prop}

\begin{pf} (1) follows from Theorem~3.1 in \cite{App} and (2)
from Proposition~6.6.1 in \cite{Far}. (3) is a straightforward
consequence of the Sobolev embedding theorem.
\end{pf}
\subsection{Examples} Now we consider different families of
measures and apply Proposition~\ref{main}. In all cases, we take
$\mu\in\CID$ so that $c_{\pi} =
e^{-\eta(\kappa_{\pi}^{\cfrac{1}{2}})}$.

\subsubsection{The case where there is a nontrivial Gaussian
component}
We say that $\mu$ has a \textit{nontrivial Gaussian
component}\vadjust{\goodbreak} if $\eta$ is such that $\sigma> 0$ in (\ref{Ls}). We
can obtain many examples of such measures by defining $\mu=
\mu_{1} * \mu_{2}$ where $\mu_{1}$ is Gaussian and $\mu_{2}$ is of
compound Poisson type or is obtained by subordination as in
Section~\ref{sec2}. We show that $\mu$ has a $C^{\infty}$-density for all
$\sigma> 0$. To prove this, we use (\ref{Sobcen}), (\ref{w1}),
(\ref{w2}) and (\ref{inteq}) and the fact that $\eta(u)
\geq\frac{1}{2}\sigma^{2}u^{2}$ for all $u \in\R$ to see that for
all $k \in\mathbb{N}$
\begin{eqnarray*}\sum_{\pi\in\widehat{G}}d_{\pi}^{2}(1 +
\kappa_{\pi})^{k}c_{\pi}^{2}
& \leq& \sum_{\pi\in
\widehat{G}}d_{\pi}^{2}(1 +
\kappa_{\pi})^{k}\exp \{-\sigma^{2}\kappa_{\pi} \}\\[-2pt]
& \leq& M\sum_{\lambda\in P \cap D}|\lambda|^{2m}(1 +
|\lambda|^{2})^{k}
\exp \{-\sigma^{2}|\lambda|^{2} \} \\[-2pt]
& \leq& K_{1}\sum_{n \in\mathbb{Z}^{r}}\|n\|^{2m}(1 +
\|n\|^{2})^{k}
\exp \{-K_{2}\|n\|^{2} \}\\[-2pt]
& = & K_{1}\sum_{j=0}^{\infty}a(j)j^{m}(1 +
j)^{k}
\exp \{-K_{2}j \}\\[-2pt]
& \leq& K_{1}\sum_{j=0}^{\infty}j^{m}\bigl(2 \sqrt{j} + 1\bigr)^{r}(1 +
j)^{k} \exp \{-K_{2}j \} < \infty,
\end{eqnarray*}
where $M, K_{1},
K_{2} > 0$, $a(j):= \#\{n \in\mathbb{Z}^{r}; \|n\|^{2} = j\}$ and
we use the fact that $a(j) \leq(2 \sqrt{j} + 1)^{r}$ for all $j
\in\mathbb{N}$.\footnote{Of course in the pure Gaussian case,
smoothness of the density is well known and can be proved using
pde techniques.}

\vspace*{-2pt}\subsubsection{Stable-type densities}\label{sec3.4.2}
Take $c_{\pi} =
\exp \{-b^{\alpha}\kappa_{\pi}^{\cfrac{\alpha}{2}} \}$
with $0 < \alpha< 2$. Again we show that the densities are
$C^{\infty}$. Indeed arguing as above
we have
\begin{eqnarray*}\sum_{\pi\in\widehat{G}}d_{\pi}^{2}(1 +
\kappa_{\pi})^{k}c_{\pi}^{2}
& \leq& M \sum_{\lambda\in P \cap
D}|\lambda|^{2m}(1 +
|\lambda|^{2})^{k}\exp \{-2b^{\alpha}|\lambda|^{\alpha}
\}\\[-2pt]
& \leq& K_{3}\sum_{j=0}^{\infty}j^{m}\bigl(2 \sqrt{j} + 1\bigr)^{r}(1 +
j)^{k}\exp \{-K_{4}j^{\cfrac{\alpha}{2}} \},
\end{eqnarray*}
where
$K_{3}, K_{4} > 0$.
To see that the series converges, it is sufficient to show that $\sum
_{n=1}^{\infty}n^{\kappa}e^{-n^{\beta}}$ converges for all $\kappa
\geq0$ where $0 < \beta< 1$. This follows by comparison with $\sum
_{n=1}^{\infty}\frac{1}{n^{2}}$ since \mbox{$\lim_{n \rightarrow\infty
}n^{\kappa+ 2}e^{-n^{\beta}} =
\lim_{x \rightarrow\infty}x^{\bfrac{\kappa+ 2}{\beta}}e^{-x} = 0$}.

\subsubsection{Relativistic Schr\"{o}dinger density}\label{sec3.4.3}

Here we have
\begin{eqnarray*}\sum_{\pi\in\widehat{G}}d_{\pi}^{2}(1 +
\kappa_{\pi})^{k}c_{\pi}^{2} & \leq& e^{2m}\sum_{\lambda\in P
\cap D}|\lambda|^{2m}\bigl|1 +
|\lambda|^{2}\bigr|^{k}e^{-2\sqrt{m^{2} + |\lambda|^{2}}}\\[-2pt]
& \leq& e^{2m}\sum_{\lambda\in P \cap D}|\lambda|^{2m}\bigl|1 +
|\lambda|^{2}\bigr|^{k}e^{-2|\lambda|} < \infty,
\end{eqnarray*}
so by the result
of Section~\ref{sec3.4.2} (with $\alpha = 1$) this case also yields a
$C^{\infty}$ density.

\subsubsection{Laplace density}
In this case, we take $c_{\pi} = (1 + \beta^{2}
\kappa_{\pi})^{-1}$. We restrict ourselves to seeking an
$L^{2}$-density. Applying (\ref{L2}), we use (\ref{w1}), (\ref{w2})
and (\ref{w3}) to obtain\vspace*{-1pt}
\[\sum_{\pi\in\widehat{G - \{\delta\}}} \frac
{d_{\pi}^2}{(1
+ \beta^2 \kappa_{\pi})^{2}} \leq\frac{N}{\beta^{2}}
\sum_{\lambda\in P \cap D-\{0\}}
\frac{|\lambda|^{2m}}{|\lambda|^4}  =
\frac{N}{\beta^{2}}\zeta(2-m),
\]
where $N > 0$.

By Sugiura's convergence result for $\zeta(s)$, we see that a
sufficient condition for convergence is $m < 2 - \frac{r}{2}$.
Hence by (\ref{rank}), $\operatorname{dim}(G) \in\{1,2,3\}$. So for
example, the Laplace distribution has a square-integrable density
on the groups $\operatorname{SO}(3), \operatorname{SU}(2)$ and $\operatorname{Sp}(1)$,
each of which has dimension $3$ and rank $1$.

\section{Pseudo-differential operator representations}\label{sec4}

In this section, $G$ is an arbitrary compact group. Let $(\mu_{t},
t \geq0)$ be a vaguely continuous convolution semigroup of
probability measures on $G$ wherein $\mu_{0} = \delta_{e}$. It
then follows that~$\mu_{t}$ is infinitely divisible for each $t
\geq0$. We let $(T_{t}, t \geq0)$ be the associated~$C_{0}$
semigroup on~$C(G)$ (\textit{Hunt semigroup}) defined by
\[
T_{t}f(\sigma) = \int_{G}f(\sigma\tau) \mu_{t}(d\tau)
\]
for all $t \geq0$. Necessary and sufficient conditions for a
densely defined linear operator to extend to the infinitesimal
generator of $(T_{t}, t \geq0)$ were found by Hunt \cite{Hu} (see
\cite{Liao} for a modern treatment) in the case of a Lie group and
generalized by Born \cite{Bo} to locally compact groups.

$(T_{t}, t \geq0)$ extends to a positivity-preserving contraction
semigroup on  $L^{2}(G): = L^{2}(G, \mathbb{C})$ and from now on we
will always work with this extended action. Our aim in this
section is to represent the semigroup and its generator as
pseudo-differential operators using Peter--Weyl theory (cf.
\cite{RT}). If $A \in M_{n}(\C)$, we define its Hilbert--Schmidt
norm by $\|A\|_{\hs}: = \operatorname{tr}(AA^{*})^{\cfrac{1}{2}}$.





The celebrated Peter--Weyl theorem asserts\vspace*{1pt} that $f \in L^{2}(G)$
has the associated Fourier series $\sum_{\pi\in
\widehat{G}}d_{\pi}\operatorname{tr}(\widehatb{f}(\pi)\pi)$ and we will make
frequent use of Plancherel's theorem in this context which tells
us that
\[
\|f\|^{2} = \sum_{\pi\in\G}d_{\pi}\|\widehatb{f}(\pi)\|^{2}_{\hs}.
\]
We also need the corresponding Parseval identity:
\[
\langle f, g \rangle= \sum_{\pi\in
\G}d_{\pi}\operatorname{tr}(\widehatb{f}(\pi)\widehatb{g}(\pi)^{*})
\]
for $f,g
\in L^{2}(G)$.\vadjust{\goodbreak}

We say that a densely defined linear operator $S$ on $L^{2}(G)$
has a (\textit{simple}) \textit{pseudo-differential operator representation}
if for each $\pi\in\widehat{G}$ there exists a~$d_{\pi} \times
d_{\pi}$ matrix $\sigma_{S}(\pi)$ such that
\[
\widehat{Sf}(\pi) = \sigma_{S}(\pi)\widehatb{f}(\pi)
\]
for all $f \in\Dom(S)$ and all $\pi\in\widehat{G}$.

We call $\sigma_{S}$ the \textit{symbol} of the operator $S$. The
word ``simple'' is included as we do not require the symbol be a
function defined on $G \times\G$ as in \cite{RT}. Indeed such a
more complicated class of symbols will be associated with
representations of more general Feller--Markov semigroups (see
\cite{Jac4} for the case where $G$ is the real numbers).

\begin{theorem} \label{pde1} For each $t \geq0, T_{t}$ is a
pseudo-differential operator with symbol $\widehata{\mu_{t}}(\pi)$ at
$\pi\in\G$.
\end{theorem}

\begin{pf} For each $\rho\in G$ let $R_{\rho}$ denote right\vspace*{1pt}
translation so that $R_{\rho}f(\sigma) = f(\sigma\rho)$ for each
$f \in L^{2}(G), \sigma\in G$. We will need the fact that
$\widehat{R_{\rho}f}(\pi) = \pi(\rho)\widehatb{f}(\pi)$ for each
$\pi\in\G$.

By Fubini's theorem and the Parseval identity,
\begin{eqnarray*}\|T_{t}f\|^{2} & = & \int_{G}\int_{G}\int
_{G}f(\sigma\tau
)\overline{f(\sigma\rho)}\mu_{t}(d\tau)\mu_{t}(d\rho)\,d\sigma\\
& = & \int_{G}\int_{G}\int_{G}R_{\tau}f(\sigma)\overline{R_{\rho
}f(\sigma)}\,d\sigma\mu_{t}(d\tau)\mu_{t}(d\rho)\\
& = & \int_{G}\int_{G}\sum_{\pi\in
\widehat{G}}d_{\pi}\operatorname{tr}(\pi(\tau)\widehatb{f}(\pi)\widehatb{f}(\pi)^{*}\pi(\rho)^{*})\mu_{t}(d\tau)\mu_{t}(d\rho).
\end{eqnarray*}

We can use Fubini's theorem to interchange summation and
integration since by the contraction property of $T_{t}$,
\[
\int_{G}\int_{G}\sum_{\pi\in
\widehat{G}}d_{\pi}\operatorname{tr}(\pi(\tau)\widehatb{f}(\pi)\widehatb{f}(\pi)^{*}\pi(\rho)^{*})\mu_{t}(d\tau)\mu_{t}(d\rho) \leq\|f\|^{2}.
\]
Hence, we have
\begin{eqnarray*}\|T_{t}f\|^{2} & = & \sum_{\pi\in\widehat
{G}}d_{\pi}\operatorname
{tr} \biggl(\int_{G}\pi(\tau)\mu_t(d\tau)\widehatb{f}(\pi)\widehatb{f}(\pi)^{*}\int_{G}\pi(\rho)^{*}\mu_t(d\rho) \biggr)\\
& = & \sum_{\pi\in
\widehat{G}}d_{\pi}\|\widehata{\mu_{t}}(\pi)\widehatb{f}(\pi)\|^{2}_{\hs},
\end{eqnarray*}
and the result follows.
\end{pf}

Let ${\cal A}$ be the infinitesimal generator of $(T_{t}, t \geq
0)$. We here use the fact that for each $t \geq0, \pi\in\G,
\widehata{\mu_{t}}(\pi) = e^{t{\cal L}_{\pi}} $ where ${\cal
L}_{\pi}$ is a $d_{\pi} \times d_{\pi}$ matrix (see\vadjust{\goodbreak} \cite{LoNg,Hey,App1}
where an explicit ``L\'{e}vy--Khintchine
type'' representation of ${\cal L}_{\pi}$ can be found when $G$ is
a Lie group).\vspace*{-3pt}

\begin{theorem} \label{pde2} ${\cal A}$ is a pseudo-differential
operator with symbol ${\cal L}_{\pi}$ at $\pi\in\G$.\vspace*{-3pt}
\end{theorem}

\begin{pf} For each $f \in\Dom({\cal A}), g \in L^{2}(G)$, we have
by Parseval's identity and Theorem~\ref{pde1}
\begin{eqnarray*}\langle{\cal A}f, g \rangle& = & \lim_{t
\rightarrow0}\sum
_{\pi\in
\G}d_{\pi}\operatorname{tr} \biggl(\frac{(\widehata{\mu_{t}}(\pi) - I_{\pi
})}{t}\widehatb{f}(\pi)\widehatb{g}(\pi)^{*} \biggr)\\[-2pt]
& = & \lim_{t \rightarrow0}\sum_{\pi\in
\G}d_{\pi}\operatorname{tr} \biggl(\frac{(e^{t{\cal L}_{\pi}} -
I_{\pi})}{t}\widehatb{f}(\pi)\widehatb{g}(\pi)^{*} \biggr).
\end{eqnarray*}

Now fix $\pi^{\prime} \in\widehat{G}$ and let $g \in{\cal
M}_{\pi^{\prime}}$ where ${\cal M}_{\pi^{\prime}}$ is the subspace
of $L^{2}(G)$ generated by mappings of the form $\sigma
\rightarrow\langle\pi^{\prime}(\sigma)u, v \rangle$ for $u, v \in
V_{\pi^{\prime}}$. It follows from the Peter--Weyl theorem that
$\widehatb{g}(\pi) = 0$ if $\pi\neq\pi^{\prime}$ and so
\begin{eqnarray*}\langle{\cal A}f, g \rangle& = & d_{\pi^{\prime
}}\lim_{t
\rightarrow0}\operatorname{tr} \biggl(\frac{(e^{t{\cal L}_{\pi^{\prime}}} -
I_{\pi^{\prime}})}{t}
\widehatb{f}(\pi^{\prime})\widehatb{g}(\pi^{\prime})^{*} \biggr)\\[-2pt]
& = & d_{\pi^{\prime}}\operatorname{tr}({\cal
L}_{\pi^{\prime}}\widehatb{f}(\pi^{\prime})\widehatb{g}(\pi
^{\prime})^{*})\\[-2pt]
& = & \sum_{\pi\in\G}d_{\pi}\operatorname{tr}({\cal
L}_{\pi}\widehatb{f}(\pi)\widehatb{g}(\pi)^{*}).
\end{eqnarray*}

The required result follows from the Parseval identity since (by
the Peter--Weyl theorem) $L^{2}(G)$ is the closure of
$\bigoplus_{\pi\in\G}{\cal M}_{\pi}$ (see, e.g., \cite{Far},
page~108).\vspace*{-2pt}
\end{pf}




For completeness, we will also give the pseudo-differential operator
representation of the resolvent $R_{\lambda}:=(\lambda I- {\cal
A})^{-1}$ for $\lambda> 0$.\vspace*{-2pt}

\begin{prop} For each $\lambda> 0, R_{\lambda}$ is a
pseudo-differential operator having
symbol $(\lambda I_\pi- {\cal L}_{\pi})^{-1}$ at $\pi\in\G$.\vspace*{-2pt}
\end{prop}

\begin{pf} First, note that $(\lambda I_\pi- {\cal L}_{\pi})^{-1}$
always exists since the eigenvalues of the matrix ${\cal L}_{\pi}$
have nonpositive real parts. We use the fact that for all\vspace*{1pt}
$\lambda
> 0, R_{\lambda} = \int_{0}^{\infty}e^{-\lambda t}T_{t}\,dt$. Then for
all $f,g \in L^{2}(G)$, by Theorem~\ref{pde1}
\[
\langle R_{\lambda}f, g \rangle= \sum_{\pi\in\G}d_{\pi}\int
_{0}^{\infty
}e^{-\lambda t}\operatorname{tr}(e^{t {\cal L}_{\pi}}\widehatb{f}(\pi
)\widehatb{g}(\pi)^{*})\,dt.
\]

The result follows from Fubini's theorem using the estimate
\[
\sum_{\pi\in\G}d_{\pi}\int_{0}^{\infty}e^{-\lambda t}|\operatorname
{tr}(e^{t {\cal L}_{\pi}}\widehatb{f}(\pi)\widehatb{g}(\pi)^{*})|\,dt
\leq\frac{1}{\lambda}\|f\|\|g\|,
\]
which is obtained by routine\vadjust{\goodbreak} computations.
\end{pf}

Now we assume that $\mu_{1} \in\CID$. It follows that $\mu_{t} \in
\CID$ for all $ t \geq0$ and that
$\widehata{\mu_{t}}(\pi) = e^{-t \eta(\kappa_{\pi}^{\cfrac
{1}{2}})} I_\pi$ for each $\pi\in\G$ for some negative definite function
$\eta$ and ${\cal A}$ has symbol whose value at $\pi\in\G$ is $-\eta
(\kappa_{\pi}^{\cfrac{1}{2}})I_{\pi}$.\vspace*{-2pt}

\begin{theorem} If $G$ is a compact Lie group and $\mu_{1} \in\CID$
then for all
$p \geq2, {\cal H}_{p}(G) \subseteq\Dom({\cal A})$ and ${\cal
A}$ is a bounded linear operator from ${\cal H}_{p}(G)$ to ${\cal
H}_{p-2}(G)$.\vspace*{-2pt}
\end{theorem}

\begin{pff*}[(Cf. \cite{RT}, Theorem~10.81, pages~571--572)] We will
make use of the fact that there exists $K > 0$ such that
$|\eta(u)| \leq K(1 + |u|^{2})$ for all $u \in\R$ (see, e.g.,
\cite{Abook}, page~31). By Theorem~\ref{pde2} and (\ref{Sob}) for each
$f \in{\cal H}_{p}(G)$
\begin{eqnarray*}\vertvertvert {\cal A}f\vertvertvert _{p-2}^{2} & = & \sum_{\pi\in\G
}d_{\pi}(1 +
\kappa_{\pi})^{p-2}\|{\cal L}_{\pi}\widehatb{f}(\pi)\|_{\hs}^{2}\\
& = & \sum_{\pi\in\G}d_{\pi}(1 + \kappa_{\pi})^{p-2}|\eta
(\kappa_{\pi}^{\cfrac{1}{2}})|^{2}\|\widehatb{f}(\pi)\|_{\hs}^{2}\\
& \leq& K\sum_{\pi\in\G}d_{\pi}(1 + \kappa_{\pi})^{p}\|\widehatb{f}(\pi)\|_{\hs}^{2}\\
& = & K\vertvertvert f\vertvertvert _{p}^{2}.
\end{eqnarray*}
In particular, it follows that
$\|{\cal A}f\|^{2} < \infty$ and so $f \in\Dom({\cal A}).$\qed\vspace*{-2pt}
\end{pff*}
\section{Transition densities for convolution semigroups}\label{sec5}

In this section, we continue to work with the Hunt semigroup
$(T_{t}, t \geq0)$ acting on the space~$L^{2}(G)$ that is
associated to the convolution semigroup of measures $(\mu_{t}, t
\geq0)$ on a compact group $G$. Let $L^{2}_{c}(G):=\{f \in
L^{2}(G), f(g \sigma g^{-1}) = f(\sigma) $  for all  $ \sigma,g
\in G\}$. It is well known that $\{\chi_{\pi}, \pi\in
\widehat{G}\}$ is a complete orthonormal basis for $L^{2}_{c}(G)$
(see, e.g., \cite{Far}, Proposition~6.5.3, page~117).\vspace*{-2pt}

\begin{prop} If $\mu_{t} \in M_{c}(G)$ for some $t \geq0$, then
$T_{t}(L^{2}_{c}(G)) \subseteq
L^{2}_{c}(G)$.\vspace*{-2pt}
\end{prop}

\begin{pf} For all $\sigma, g \in G, f \in L^{2}_{c}(G)$, we have
\begin{eqnarray*}T_{t}f(g \sigma g^{-1}) & = & \int_{G}f(g \sigma g^{-1}
\tau)
\mu_{t}(d\tau)\\[-2pt]
& = & \int_{G}f(g \sigma g^{-1} \tau) \mu_{t}(d g^{-1}\tau g)\\[-2pt]
& = & \int_{G}f(g \sigma\tau g^{-1}) \mu_{t}(d\tau)\\[-2pt]
& = & T_{t}f(\sigma).\vspace*{-2pt}
\end{eqnarray*}
\upqed
\end{pf}

Now suppose that $\mu_{t} \in M_{cs}(G)$ for all $t \geq0$. This
implies in particular that $(T_{t}, t \geq0)$ is self-adjoint in
$L^{2}(G)$ and hence in $L^{2}_{c}(G)$ (see \cite{App2,Kun}). By Proposition~\ref{non-zero}, we also have that
there exists $\alpha_{\pi} \leq0$ for each $\pi\in\widehat{G}$
such that $\widehata{\mu_{t}}(\pi) = e^{t \alpha_{\pi}}I_\pi$.

\begin{theorem} \label{comp} If $\mu_{t} \in M_{cs}(G)$ for all
$t \geq0$, then $\{\chi_{\pi}, \pi\in\widehat{G}\}$ is a
complete set of eigenvectors for the action of $T_{t}$ on
$L^{2}_{c}(G)$ and
\begin{equation} \label{ev}
T_{t}\chi_{\pi} = e^{t \alpha_{\pi}}\chi_{\pi}
\end{equation}
for all $\pi\in\widehat{G}, t \geq0$.
\end{theorem}

\begin{pf} For all $\sigma\in G, t \geq0$,
\begin{eqnarray*}T_{t}
\chi_{\pi} (\sigma) & = & \int_{G} \chi_{\pi}(\sigma\tau)
\mu_{t}(d\tau)\\
& = & \int_{G} \operatorname{tr}(\pi(\sigma) \pi(\tau)) \mu_{t}(d\tau)\\
& = & \operatorname{tr} \biggl( \pi(\sigma) \int_{G} \pi(\tau) \mu
_{t}(d\tau) \biggr)\\
& = & e^{t \alpha_{\pi}}\chi_{\pi}(\sigma).
\end{eqnarray*}
\upqed
\end{pf}

It is shown in \cite{App2} that $T_{t}$ is trace-class for $t > 0$
if and only if $\mu_{t}$ has a square-integrable density. In this
case, we have
\begin{equation} \label{tra}
\operatorname{tr}(T_{t}) = \sum_{\pi\in\G}e^{t \alpha_{\pi}}
\end{equation}
for $t > 0$, where $\operatorname{tr}$ denote the trace in the Hilbert space
$L^{2}_{c}(G)$.

From now on, we assume that for $t > 0, \mu_{t} \in{\cal
M}_{cs}(G)$ has a density $k_{t} \in L^{2}_{c}(G, \R)$. We define
the \textit{transition density} $h_{t} \in L^{2}(G \times G, \R)$ by
\[
h_{t}(\sigma, \tau): = k_{t}(\sigma^{-1}\tau)
\]
for each $t > 0, \sigma,\tau\in G$. Indeed $h_{t}$ is precisely
the transition probability density at time $t$ of a $G$-valued
L\'{e}vy process whose law at time $t$ is $k_{t}$.

Note that for each $\sigma\in G, (t, \rho) \rightarrow
h_{t}(\sigma, \rho)$ satisfies the backward equation (in the
distributional sense)
\[
\frac{\partial h_{t}}{\partial t}(\sigma, \rho) = {\cal
A}h_{t}(\sigma, \rho),
\]
with $h_{t}(\sigma, \rho) \rightarrow\delta_{\sigma}(\rho)$ as $t
\rightarrow0$.
For example, if $(\mu_{t}, t \geq0)$ is the Brownian motion semigroup
which is characterized
by $\widehata{\mu}_{t}(\pi) = e^{-\afrac{t}{2}\kappa_{\pi}}I_{\pi}$
for each $\pi\in\widehat{G}$, then $h_{t}$ is the well-known
heat kernel and for this reason we will sometimes refer to our
more general $h_{t}$ as the \textit{transition kernel}.

\begin{theorem} For each $t > 0$:
\begin{longlist}[(2)]
\item[(1)]
\begin{equation} \label{td1} \int_{G}h_{t}(g^{-1}\sigma,
\rho g^{-1})\,dg = \sum_{\pi\in\G}e^{t
\alpha_{\pi}}\overline{\chi_\pi(\sigma)}\chi_\pi(\rho)
\end{equation}
for all $\sigma, \rho\in G$.
\item[(2)]
\begin{equation} \label{td2} \operatorname{tr}(T_{t}) =
\int_{G}\int_{G}k_{t}(\rho^{-1}g \rho g^{-1})\,dg\,d\rho.
\end{equation}
\end{longlist}
\end{theorem}

\begin{pf}
\begin{longlist}[(2)]
\item[(1)] By (\ref{dens})
\begin{equation} \label{katye}
h_{t}(\sigma, \rho) = k_{t}(\sigma^{-1}\rho) = \sum_{\pi\in
\G}d_{\pi}e^{t \alpha_{\pi}} \chi_{\pi}(\sigma^{-1}\rho),
\end{equation}
and so
\begin{eqnarray*}\int_{G}h_{t}(g^{-1}\sigma, \rho g^{-1})\,dg & = &
\sum_{\pi
\in\G}d_{\pi}e^{t \alpha_{\pi}}\int_{G}\chi_{\pi}(\sigma
^{-1}g\rho g^{-1})\,dg\\
& = & \sum_{\pi\in
\G}e^{t\alpha_{\pi}}\overline{\chi_{\pi}(\sigma)}\chi_\pi(\rho),
\end{eqnarray*}
by Proposition~6.5.2 in \cite{Far} (page~116). The interchange
of integral and sum is justified by Fubini's theorem since
\[
\sum_{\pi\in\G}d_{\pi}e^{t
\alpha_{\pi}}\int_{G}|\chi_{\pi}(\sigma^{-1}g\rho g^{-1})|\,dg \leq
\sum_{\pi\in\G}d_{\pi}^{2}e^{t \alpha_{\pi}} < \infty
\]
%
by (\ref{L2}) since $k_{t} \in L^{2}(G)$ for each $t > 0$. Here, we
have used
the crude estimate $\sup_{\sigma\in
G}|\chi_{\pi}(\sigma)| \leq d_{\pi}$.
\item[(2)] Put $\rho= \sigma$ in (\ref{td1}) and then integrate both
sides with respect to $\sigma$ using the fact that
$\int_{G}|\chi_{\pi}(\sigma)|^{2}\,d\sigma= 1$. The result then
follows from (\ref{tra}). Note that the interchange of integral and
summation is justified by Fubini's theorem using a similar argument to
that presented in (1).\qed
\end{longlist}
\noqed
\end{pf}
\begin{cor} If $k_{t}$ is continuous for each $t > 0$,
\begin{equation} \label{katy1}
k_{t}(e) = \sum_{\pi\in\G}d_{\pi}^{2}e^{t\alpha_{\pi}}.
\end{equation}
\end{cor}

\begin{pf} Put $\sigma= \rho$ in (\ref{katye}) [or argue
directly from\vadjust{\goodbreak} (\ref{dens})].
\end{pf}

We now work in the Hilbert space $L^{2}(G)$ and we use $\Tr$
to denote the trace in this Hilbert space. By the Peter--Weyl
theorem, $\{d_{\pi}^{\cfrac{1}{2}}\pi_{ij}, 1 \leq i,j \leq d_{\pi},
\pi\in\G\}$ is a complete orthonormal basis for $L^{2}(G)$ where
$\pi_{ij}$ denotes the coordinate function $\pi_{ij}(\sigma) :=
\pi(\sigma)_{ij}$, for each $\sigma\in G$.

The following two results are well known for the heat kernel (see, e.g.,
Chapter~12 of \cite{Feg}). Here we extend them to more general Hunt
semigroups.\

\begin{theorem} \label{late}
\begin{longlist}[(2)]
\item[(1)] For each $t \geq0$, the set $\{\pi_{ij}, 1 \leq i,j \leq d_{\pi
}, \pi\in\G\}$ is a complete orthogonal set of eigenvectors for $T_{t}$
and
\begin{equation} \label{late1}
T_{t}\pi_{ij} = e^{t\alpha_{\pi}}\pi_{ij}
\end{equation}
for each $1 \leq i,j \leq d_{\pi}, \pi\in\G.$
\item[(2)] If $k_{t}$ is continuous for each $t > 0$,
\begin{equation} \label{late2}
k_{t}(e) = \Tr(T_{t}).
\end{equation}
\end{longlist}
\end{theorem}

\begin{pf}
\begin{longlist}[(2)]
\item[(1)] For all $\sigma\in G$
\begin{eqnarray*}T_{t}\pi_{ij}(\sigma) & = & \int_{G}\pi
_{ij}(\sigma\tau
)\mu_{t}(d\tau)\\
& = & \sum_{k=1}^{d_{\pi}}\pi_{ik}(\sigma)\int_{G}\pi_{kj}(\tau
)\mu_{t}(d\tau)\\
& = & \sum_{k=1}^{d_{\pi}}\pi_{ik}(\sigma)\widehata{\mu_{t}}(\pi
)_{kj}\\
& = & e^{t\alpha_{\pi}}\pi_{ij}(\sigma).
\end{eqnarray*}
\item[(2)] Since each eigenvalue $e^{t\alpha_{\pi}} $ has multiplicity
$d_{\pi}^{2}$ on the closed subspace of $L^{2}(G)$ spanned by
$\{\pi_{ij}, 1 \leq i,j \leq d_{\pi}\}$ it is clear that
\[
\Tr(T_{t}) = \sum_{\pi\in\G}d_{\pi}^{2}e^{t\alpha_{\pi}},
\]
and the result follows from (\ref{katy1}).\qed
\end{longlist}
\noqed
\end{pf}

The next result can be deduced directly from (\ref{dens}). We give
a direct proof to make the paper more self-contained. Note that
results of this type are well known for Markov processes taking
values in   compact metric spaces and having suitably
square-integrable transition probabilities (see Theorem~6.4 in
\cite{Get}).

\begin{theorem} \label{eigr}
If $k_{t}$ is continuous for each $t > 0$,
\item
\begin{equation} \label{eigr1}
h_{t}(\sigma, \rho) = \sum_{\pi\in\G}\sum_{i,j = 1}^{d_{\pi
}}d_{\pi}e^{t\alpha_{\pi}} \overline{\pi_{ij}(\sigma)}\pi
_{ij}(\rho)
\end{equation}
for all $\sigma, \rho\in G$.
\end{theorem}

\begin{pf}
For each $\sigma\in G$ let $L_{\sigma}$ denote left translation so
that
$L_{\sigma}f (\rho) = f(\sigma^{-1}\rho)$ for each
$f \in L^{2}(G), \rho\in G$.
By Fourier expansion and using (\ref{late1}),
\begin{eqnarray*}L_{\sigma}k_{t} & = & \sum_{\pi\in\G}\sum_{i,j =
1}^{d_{\pi
}}d_{\pi}
\langle L_{\sigma}k_{t}, \pi_{ij} \rangle\pi_{ij} \\
& = & \sum_{\pi\in\G}\sum_{i,j = 1}^{d_{\pi}}d_{\pi} \biggl(\int
_{G}k_{t}(\sigma^{-1}\tau)\overline{\pi_{ij}(\tau)}\,d\tau
\biggr)\pi_{ij}\\
& = & \sum_{\pi\in\G}\sum_{i,j = 1}^{d_{\pi}}d_{\pi} \biggl(\int
_{G}k_{t}(\tau)\overline{\pi_{ij}(\sigma\tau)}\,d\tau \biggr)\pi
_{ij}\\
& = & \sum_{\pi\in\G}\sum_{i,j = 1}^{d_{\pi}}d_{\pi
}T_{t}\overline{\pi_{ij}(\sigma)}\pi_{ij}\\
& = & \sum_{\pi\in\G}\sum_{i,j=1}^{d_{\pi}}d_{\pi}e^{t\alpha
_{\pi}}\overline{\pi_{ij}(\sigma)}\pi_{ij}.
\end{eqnarray*}

Since $\sup_{\sigma,\rho\in G} |\sum_{i,j=1}^{d_{\pi
}}\overline{\pi_{ij}(\sigma)}\pi_{ij}(\rho) | = \sup
_{\sigma,\rho\in G}|\operatorname{tr}(\pi(\sigma^{-1}\rho)| \leq d_{\pi
}$, we
deduce uniform convergence of the series from (\ref{L2}) and so for
all $\rho\in G$
\[
h_{t}(\sigma, \rho) = L_{\sigma}k_{t}(\rho) = \sum_{\pi\in\G
}\sum_{i,j = 1}^{d_{\pi}}d_{\pi}e^{t\alpha_{\pi}} \overline{\pi
_{ij}(\sigma)}\pi_{ij}(\rho).
\]
\upqed
\end{pf}

\begin{cor} If $k_{t}$ is continuous for each $t > 0$,
\[
k_{t}(e) \geq\int_{G}\int_{G}k_{t}(\rho^{-1}g \rho
g^{-1})\,dg\,d\rho
\]
with equality if and only if $G$ is Abelian.
\end{cor}

\begin{pf} The inequality follows from (\ref{tra}), (\ref{td2}),
(\ref{katy1}) and (\ref{late2}). If $G$ is Abelian, then equality
is obvious. If $G$ is non-Abelian, we must have $d_{\pi} > 1$ for
at least one $\pi\in\G$ and then it is clear that strict inequality
holds.
\end{pf}


\section{Small time asymptotics of densities}\label{sec6}

Assume that $G$ is a compact semisimple Lie group. We would like
to obtain an asymptotic expansion for $k_{t}(e)$ as $t \rightarrow
0$. We assume that $\mu_{t} \in\CID$ for each $t \geq0$. In this
case, if~$k_{t}$ is continuous for $t > 0$, we may follow the
arguments on page~106 of~\cite{Feg} to obtain
\begin{eqnarray} \label{weitra} k_{t}(e) &=& \Tr(T_{t})
 =  \sum_{\lambda\in P \cap D}d_{\lambda+
\rho}^{2}\exp\bigl\{
-t\eta\bigl((|\lambda+ \rho|^{2} -
|\rho|^{2})^{\cfrac{1}{2}}\bigr)\bigr\} \nonumber\\
 &=&  \sum_{\lambda\in P \cap
D}d_{\lambda}^{2}\exp\bigl\{-t\eta\bigl((|\lambda|^{2} -
|\rho|^{2})^{\cfrac{1}{2}}\bigr)\bigr\} \\
& = & \frac{1}{|W|} \sum_{\lambda\in
P}d_{\lambda}^{2}\exp\bigl\{-t\eta\bigl((|\lambda|^{2} -
|\rho|^{2})^{\cfrac{1}{2}}\bigr)\bigr\}, \nonumber
\end{eqnarray}
where
$d_{\lambda}$ denotes the dimension of the representation space
with highest weight $\lambda$ and $|W|$ is the order of the Weyl
group of $G$.


When $\eta(u) = \frac{|u|^{2}}{2}$, $k_{t}$ is the density
generating the heat kernel and it is known that as $t \rightarrow
0$,
\begin{equation} \label{hkas}
k_{t}(e) \sim C t^{-\bfrac{\operatorname{dim}(G)}{2}}e^{t|\rho|^{2}}
\end{equation}
(see, e.g., \cite{Feg}, page~109), where $C > 0$.

We will examine the case where $\eta$ is the characteristic
exponent of a~symmetric Cauchy distribution so that $\eta(u) =
\sigma|u|$ for all $u \in\R$, where $\sigma> 0$.

\begin{exampleb} $G = \Pi^{d}$ where $\Pi:= \R/\mathbb{Z}$.
In this case $\widehat{G} = \mathbb{Z}^{d}$, each $d_{\pi} = 1$
and for each $\pi= n \in\mathbb{Z}^{d}, \kappa_{\pi} =
4\pi^{2}|n|^{2}$ where $n^{2} = n_{1}^{2} + \cdots+ n_{d}^{2}$
for $n = (n_{1}, \ldots, n_{d})$. The equation (\ref{weitra}) then
takes the form
\[
k_{t}(e) = \sum_{n \in\mathbb{Z}^{d}}e^{-2\pi t\sigma|n|}.
\]

When $d = 1$, we easily calculate
\begin{eqnarray*}k_{t}(e) & = & 1 + 2\sum_{n=1}^{\infty}e^{-2\pi
t\sigma n}\\
& = & \coth(\pi\sigma t) \sim\frac{1}{\pi\sigma
t} \qquad  \mbox{as } t \rightarrow0.
\end{eqnarray*}

When $d > 1$, we apply the Poisson summation formula to obtain
\begin{eqnarray*}k_{t}(e) & = &
\frac{\Gamma(\bfrac{d+1}{2})}{\pi^{\bfrac{d+1}{2}}}\sum_{m \in
\mathbb{Z}^{d}}\frac{\sigma t}{(\sigma^{2}t^{2} +
|m|^{2})^\bfrac{d+1}{2}}\\
& \sim&
\frac{\Gamma(\bfrac{d+1}{2})}{\sigma^{d}\pi^{\bfrac
{d+1}{2}}}\frac{1}{t^{d}} \qquad  \mbox{as }
t \rightarrow0.
\end{eqnarray*}
\end{exampleb}

\begin{exampleb}\label{exa2} $G = \operatorname{SU}(2)$. In this case, $\G\cong
\mathbb{Z}_{+}$ with $d_{n} = n+1$ and for each $n \in
\mathbb{Z}_{+}, \kappa_{n} = n(n+2)$. Hence,
\begin{eqnarray*}k_{t}(e) & = & \sum_{n=0}^{\infty
}(n+1)^{2}e^{-t\sigma\sqrt
{n(n+2)}}\\
& = & \sum_{n=0}^{\infty}(n+1)^{2}e^{-t\sigma\sqrt{(n+1)^{2} - 1}}\\
& = & \sum_{m = 1}^{\infty}m^{2}e^{-t\sigma\sqrt{m^{2} - 1}}.
\end{eqnarray*}

From this, we get the easy estimate
\begin{equation} \label{basicest}
e^{-t\sigma}\sum_{m = 1}^{\infty}m^{2}e^{-t\sigma m} \leq k_{t}(e)
\leq e^{t\sigma}\sum_{m = 1}^{\infty}m^{2}e^{-t\sigma m}.
\end{equation}

Now define $ g(t): = \sum_{m = 1}^{\infty}e^{-t\sigma m}$ for $t
\in(0, \infty)$. The function $g$ is $C^{\infty}$ and we have
\begin{eqnarray*}\sum_{m=1}^{\infty}m^{2}e^{-t\sigma m} & = &
\frac{1}{\sigma^{2}}\frac{d^{2}}{dt^{2}}g(t)\\
& = &
\frac{1}{\sigma^{2}}\frac{d^{2}}{dt^{2}} \biggl(\frac
{e^{-t\sigma}}{1
- e^{-t
\sigma}} \biggr)\\
& = & \frac{e^{-t\sigma}}{(1 - e^{-t
\sigma})^{2}}\coth \biggl(\frac{\sigma t}{2} \biggr),
\end{eqnarray*}
and hence we conclude that
\[
k_{t}(e) \sim\frac{2}{\sigma^{3}t^{3}} \qquad  \mbox{as } t \rightarrow0.
\]

This should be compared with the usual heat kernel where the
following very precise asymptotic expansion is known (see
\cite{Feg2}, Proposition~2.3, page~662):
\[
k_{t}(e) \sim32\sqrt{2}\pi^{2}(4\pi\sigma t)^{-\cfrac
{3}{2}}e^{\bfrac{\sigma t}{8}},
\]
so the leading term has the slower decay
\[
k_{t}(e) \sim32\sqrt{2}\pi^{2}(4\pi\sigma)^{-\cfrac{3}{2}}
t^{-\cfrac{3}{2}}.
\]

It is not difficult to verify that the relativistic
Schr\"{o}dinger semigroup on $\operatorname{SU}(2)$ with mass parameter $m$, for
which $\eta(u) = (m^{2} + u^{2})^{\cfrac{1}{2}}- m$ for $u \in\R$,
has exactly the same short time asymptotics as the Cauchy
semigroup.\looseness=1
\end{exampleb}

\begin{exampleb} $G=S0(3)$. Here we have $\G\cong\mathbb{Z}_{+}$
with $d_{n} = 2n+1$ and for each $n \in\mathbb{Z}_{+}, \kappa_{n}
= n(n+1)$. So we obtain
\begin{eqnarray*}k_{t}(e) & = &
\sum_{n=0}^{\infty}(2n+1)^{2}e^{-t\sigma\sqrt{n(n+1)}}\\
& = &
\sum_{m\ \mathrm{odd}}m^{2}e^{-t\afrac{\sigma}{2}\sqrt{m^{2}-1}} \qquad  \mbox{as } t
\rightarrow0.
\end{eqnarray*}

Using the results obtained in Example~\ref{exa2}, we find that
\[
k_{t}(e) \sim\frac{8}{\sigma^{3}t^{3}} \qquad  \mbox{as } t \rightarrow0.
\]

Based on these calculations, we conjecture that $k_{t}(e) \sim
Ct^{-d}$ for the Cauchy process on an arbitrary compact semisimple
Lie group of dimension $d$. It is also tempting to further
conjecture that if $k_{t}$ is associated to an arbitrary
$\alpha$-stable process so that $\eta(u) = |u|^{\alpha}$ where $0
< \alpha\leq2$, then $k_{t}(e) \sim Ct^{-\cfrac{d}{\alpha}}$ as
$t \rightarrow0$. This is consistent with the known behavior of
densities of symmetric stable processes in Euclidean space (see,
e.g., \cite{BG2}) where it essentially follows by scaling
arguments. We remind the reader that this technology is not
available on compact groups (see, e.g., Theorem~2.2. in \cite{App3},
page~117).

Now let $N(\lambda)$ denote the number of eigenvalues of $-{\cal
A}$ that do not exceed~$\lambda$ and note that for all $t > 0$,
\[
\Tr(T_{t}) = \int_{0}^{\infty}e^{-t
\lambda}\,dN(\lambda).
\]

If the above conjecture holds then by Karamata's Tauberian theorem
we have
\[
N(\lambda) \sim\frac{C
\lambda^{\cfrac{d}{\alpha}}}{\Gamma (1 +
\frac{d}{\alpha} )} \qquad  \mbox{as } \lambda \rightarrow\infty
\]
(cf.
Theorem~2.3 in \cite{BG1}).

So far we know that this eigenvalue asymptotics is valid when
$\alpha= 2$, and when $\alpha= 1$ with $G = \Pi^{d}, G = \operatorname{SU}(2)$
and $G = \operatorname{SO}(3)$.\vspace*{3pt}
\end{exampleb}

\section*{Acknowledgments}  I would like to thank Natesh Pillai for
many useful discussions and Ming Liao and Michael Ruzhansky for
helpful comments. Both Ren\'{e} Schilling and Rodrigo Ba\~{n}uelos
provided very valuable remarks after I presented a talk based on
part of this paper at the 2010 L\'{e}vy processes conference in
Dresden. Last, but not least, I am grateful to the referee for his
careful reading and a number of helpful observations.\looseness=1

%

\printaddresses

\end{document}